\documentclass[preprint,12pt]{elsarticle}
\usepackage{amssymb}
\usepackage{mathrsfs}
\usepackage{epic,eepic,epsf,epsfig}
\usepackage{graphicx}
\usepackage{amsfonts,srcltx,mathrsfs}
\usepackage{amsmath}
\usepackage{amsbsy}
\usepackage{amsfonts}%,srcltx}
\usepackage{color}
\usepackage{array}

\usepackage{amssymb}
\usepackage{amsthm}
\newtheorem{theorem}{Theorem}
\newtheorem{lem}{Lemma}
\newtheorem{cor}{Corollary}
\newdefinition{rmk}{Remark}
\newdefinition{defi}{Definition}
\newproof{pf}{Proof}

\def\f{\noindent}
\def\Cay{\hbox{\rm Cay}}
\def\demo{{\bf Proof.}\hskip10pt}

\journal{}

\begin{document}

\begin{frontmatter}

\title{Equal relation between the extra connectivity and pessimistic
diagnosability for some regular graphs}

\author[bjtu]{Mei-Mei Gu}
\ead{12121620@bjtu.edu.cn,}

\author[bjtu]{Rong-Xia Hao\corref{cor1}}
\ead{rxhao@bjtu.edu.cn,}

\author[ustc]{Jun-Ming Xu}
\ead{xujm@ustc.edu.cn}

\author[bjtu]{Yan-Quan Feng}
\ead{yqfeng@bjtu.edu.cn,}

\address[bjtu]{Department of Mathematics, Beijing Jiaotong University, Beijing, 100044, China}
\address[ustc]{School of Mathematical Sciences, University of Science and Technology of China, Hefei, Anhui, 230026, China}

%\author{ }
%\ead{}
%
%\author{ \corref{cor1}}
%\ead{}
%
%
%\address{Department of Mathematics, Beijing Jiaotong University
%, 100044, China}
%
\cortext[cor1]{Corresponding author}

\begin{abstract}
Extra connectivity and the pessimistic diagnosis are two crucial subjects
for a multiprocessor system's ability to tolerate
and diagnose faulty processor. The pessimistic diagnosis strategy is a classic strategy based
on the PMC model in which isolates all faulty vertices
within a set containing at most one fault-free vertex.
In this paper, the result that the pessimistic
diagnosability $t_p(G)$ equals the extra connectivity $\kappa_{1}(G)$ of a regular graph $G$
under some conditions are shown. Furthermore, the following new results are gotten: the pessimistic diagnosability $t_p(S_n^2)=4n-9$ for split-star networks $S_n^2$; $t_p(\Gamma_n)=2n-4$ for Cayley graphs generated by transposition trees $\Gamma_n$; $t_p(\Gamma_{n}(\Delta))=4n-11$ for Cayley graph generated by the $2$-tree $\Gamma_{n}(\Delta)$; $t_{p}(BP_n)=2n-2$ for the burnt pancake networks $BP_n$. As corollaries, the known results about the extra connectivity and the pessimistic diagnosability of many famous networks including the alternating group graphs; the alternating group networks; BC networks; the $k$-ary $n$-cube networks etc. are obtained directly.

\end{abstract}


\end{frontmatter}

\section{Introduction}

It is well known that a topological structure of an interconnection network can be modeled by a loopless undirected graph $G=(V,E)$, where vertices in $V$ represent the processors and the edges in $E$ represent the communication links. In this paper, we use graphs and networks interchangeably. The {\it connectivity} $\kappa(G)$ of a connected graph $G$ is the minimum number of vertices removed to
get the graph disconnected or trivial. In a multiprocessor system, some processors may fail, connectivity is used to determine the reliability and fault tolerance of a network. However, connectivity is not suitable for large-scale processing systems because it is almost impossible for all processors adjacent to, or all links incident to, the same processors to fail simultaneously. To
compensate for this shortcoming, it seems reasonable to
generalize the notion of classical connectivity by imposing
some conditions or restrictions on the components
of $G$ when we delete the set of faulty processors.
J. F\'abrega and M.A. Fiol \cite{F¡®abrega} introduced the {\it extra connectivity} of interconnection networks as follows.

\begin{defi}
A vertex set $S\subseteq V(G)$ is called to be an {\it$h$-extra vertex cut} if $G-S$ is disconnected and every component of $G-S$ has at least $h+1$ vertices. The {\it$h$-extra connectivity} of $G$, denoted by $\kappa_{h}(G)$, is defined as the cardinality of a minimum $h$-extra vertex cut, if exists.
\end{defi}

It is obvious that $\kappa_{0}(G)=\kappa(G)$ for any graph $G$ that is not a complete graph. The $1$-extra connectivity is usually called extra connectivity.
The problem of determining the extra connectivity of numerous networks has received a great deal of attention in recent years. Interested readers may refer to \cite{B99,G,G15,Hsieh,LX15,LZX} or others for further details.

The diagnosis of a system is
the process of appraising the faulty processors.
A number of models
have been proposed for diagnosing faulty processors in
a network. Preparata et al.~\cite{P}
first introduced a graph theoretical model, the so-called {\it PMC model} (i.e., Preparata, Metze and Chien¡¯s model), for system level diagnosis in multiprocessor systems.
The pessimistic diagnosis strategy proposed by Kavianpour and Friedman~\cite{KK} is a classic diagnostic model based on the PMC model. In this strategy, all faulty processors to be isolated within a set having at most one fault-free processor.

\begin{defi}
A system is {\it $t/t$-diagnosable}, provided the number of faulty processors is bounded by $t$, all faulty processors can be isolated
within a set of size at most $t$ with at most one fault-free vertex mistaken as a faulty one.
The {\it pessimistic diagnosability} of a system $G$, denoted by $t_{p}(G)$, is the maximal number of faulty processors so that the system $G$ is $t/t$-diagnosable.
\end{defi}

The pessimistic diagnosability of many interconnection networks has been explored. For example, see~\cite{F,FX,KK,T1,T13,WZF} etc.

Based on the important of the extra connectivity and
the pessimistic diagnosability and motivated by the recent researches
on the extra connectivity and pessimistic diagnosability
of some graphs, including some famous networks, our object is to propose
the relationship between extra connectivity and pessimistic
diagnosability of regular graphs with some
given conditions. In this paper, the result that the pessimistic
diagnosability $t_p(G)$ equals the extra connectivity $\kappa_{1}(G)$ of a regular graph $G$
under some conditions are shown. Furthermore, the following new results are gotten: the pessimistic diagnosability $t_p(S_n^2)=4n-9$ for split-star networks $S_n^2$;
$t_p(\Gamma_n)=2n-4$ for Cayley graphs generated by transposition trees $\Gamma_n$; $t_p(\Gamma_{n}(\Delta))=4n-11$ for
Cayley graph generated by the $2$-tree $\Gamma_{n}(\Delta)$; $t_{p}(BP_n)=2n-2$ for the burnt pancake networks $BP_n$. As corollaries, the known results about the extra connectivity
and the pessimistic diagnosability of many famous networks including the alternating group graphs, the alternating group networks, BC networks and the $k$-ary $n$-cube networks etc. are obtained directly.

%and determine that the pessimistic
%diagnosability under the PMC model equals the extra connectivity under some conditions.
%As applications of the main result, the pessimistic
%diagnosability and the extra connectivity of many famous networks,
%such as the alternating group graph $AG_n$, the alternating group network $AN_n$,
%the $k$-ary $n$-cube networks $Q_n^k$, the BC networks $X_n$, the split-star networks $S_{n}^{2}$,
%the Cayley graphs generated by transposition trees $\Gamma_n$,
%the Cayley graphs generated by $2$-trees $\Gamma_n(\Delta)$, the burnt pancake networks $BP_n$, are obtained directly.

The remainder of this paper is organized
as follows. Section $2$ introduces necessary definitions and properties of some graphs.
In Section 3, we determines the equal relationship between extra connectivity and pessimistic
diagnosability of regular graphs with some
given conditions. In Section 4, we concentrates on the applications to
some famous networks. The pessimistic
diagnosability and the extra connectivity of many famous networks,
such as the alternating group graph $AG_n$, the alternating group network $AN_n$,
the $k$-ary $n$-cube networks $Q_n^k$, the BC networks $X_n$, the split-star networks $S_{n}^{2}$,
the Cayley graphs generated by transposition trees $\Gamma_n$,
the Cayley graphs generated by $2$-trees $\Gamma_n(\Delta)$ and the burnt pancake networks $BP_n$ are obtained directly.
Finally, our conclusions are given in Section 5.

\section{Preliminaries}
\f In this section, we give some terminologies and notations of combinatorial network theory.
For notations not defined here, the reader is referred to~\cite{B}.

We use a graph, denoted by $G=(V(G),E(G))$, to represent an interconnection network, where $V(G)$ is the vertex set of $G$; $E(G)$ is the edge set of $G$.
For a vertex $u\in V(G)$, let $N_{G}(u)$ (or $N(u)$ if there is no ambiguity) denote
a set of vertices in $G$ adjacent to $u$.
For a vertex set $U \subseteq V(G)$, let $N_{G}(U)=\bigcup\limits_{v\in U}N_{G}(v)-U$ and $G[U]$ be
the subgraph of $G$ induced by $U$. If $|N_{G}(u)|=k$ for any vertex in
$G$, then $G$ is {\it $k$-regular}.
For any two vertices $u$ and $v$ in $G$, let
$cn(G;u,v)$ denote the number of vertices who are the neighbors
of both $u$ and $v$, that is, $cn(G;u,v)=|N_{G}(u)\cap N_{G}(v)|$.
Let $cn(G)=\max \{cn(G;u,v): u,v\in V(G)\}$, $l(G)=\max \{cn(G;u,v): (u,v)\in E(G)\}$.
Let $|V(G)|$ be the size of vertex set and $|E(G)|$ be the size of edge set.
Throughout this paper, all graphs are finite, undirected without loops.

Let $[n]=\{1,2,\ldots,n\}$ and $\langle n\rangle=\{-1,-2,\ldots,-n,1,2,\ldots,n\}$.
For a finite group $A$ and a subset $S$ of $A$ such that $1\notin S$ and $S=S^{-1}$ (where $1$ is the identity element of $A$),
the {\it Cayley graph} $\Cay(A;S)$ on $A$ with respect to $S$ is defined to have vertex set $A$ and edge set $\{(g,gs)|g\in A, s\in S\}$. A Cayley graph is $|S|$-regular, and is connected if and only if $S$
generates $\Gamma$. Moreover, A Cayley graph is $|S|$-connected if
$S$ is a minimal generating set of $\Gamma$.

\subsection{The alternating group graphs}

Jwo et al.~\cite{jhd93} introduced the alternating group graph as an interconnection network topology for computing systems.

\begin{defi}
Let $A_n$ be the alternating group of degree $n$
with $n\geq 3$. Set $S=\{(1\ 2\ i), (1\ i\ 2)\ |\ 3\leq i\leq n\}$.
The {\it alternating group graph}, denoted by $AG_n$, is defined as
the Cayley graph $AG_n=\Cay(A_n,S)$.
\end{defi}

It is clear that $AG_3$ is a triangle, $AG_n$ is a $(2n-4)$-connected and $(2n-4)$-regular graph with
$n!/2$ vertices. Each $AG_n$ contains $n$ sub-alternating group graphs $AG_n^{0},AG_n^{1},\ldots,AG_n^{n-1}$.
For each $i\in[n]$, $AG_n^{i}$ is isomorphic to $AG_{n-1}$.
For each vertex $v\in AG_n^{i}$, $v$ has exactly two neighbors that are not contained in $AG_n^{i}$, which are called the extra neighbors of $v$.

\begin{lem}{\rm(\cite{H1})}\label{AGn1}
The extra neighbors of every vertex of $AG_{n}$ are in different subgraphs $AG_{n}^{i}$ for $n\geq4$.
For any two different vertices $u,v$, $cn(AG_{n}:u,v)=1$ if $u$ and $v$ are adjacent; otherwise, $cn(AG_{n}:u,v)\leq2$.
\end{lem}

\begin{lem}{\rm(\cite{T1})}\label{AGn3}
Let $AG_n$ be the $n$-dimensional alternating group graph for $n\geq 4$.
If $U$ is a subset of $V(AG_n)$ and $2\leq |U|\leq 8n-25$, then $|N_{AG_n}(U)|\geq 4n-11$.
\end{lem}

\begin{lem}{\rm(\cite{H1})}\label{AGn4}
Let $F$ be a vertex-cut of $AG_n$ for $n\geq 5$. If $|F|\leq 4n-11$, then $AG_n-F$ satisfies one of the following conditions:
\begin{enumerate}
\item [{\rm (1)}] $AG_n-F$ has two components, one of which is a trivial component.
\item [{\rm (2)}] $AG_n-F$ has two components, one of which is an edge.
Moreover, if $|F|=4n-11$, $F$ is formed by the neighbor of the edge.
\end{enumerate}
\end{lem}

\subsection{The alternating group networks}

The alternating group network $AN_n$
was first proposed by Y. Ji~\cite{J} to improve upon the alternating
group graph $AG_n$, studied by Jwo and others~\cite{jhd93}.

\begin{defi}{\rm(\cite{J})}
Let $A_n$ be an alternating group of degree $n\geq 3$ and let
$S=\{(1\ 2\ 3), (1\ 3\ 2),$ $(1\ 2)(3\ i)$ $|\ 4\leq i\leq n\}$. The
{\it alternating group network}, denoted by $AN_n$, is defined as
the Cayley graph $\Cay(A_n,S)$.
\end{defi}

By the definition, we can get some properties about $AN_n$~\cite{J}. $AN_n$ is a regular graph with
$n!/2$ vertices and $n!(n-1)/4$ edges. $AN_3$ is a triangle.  $AN_4$ contains four copies of $AN_{3}$.
$AN_n$ contains $n$ copies of $AN_{n-1}$, say $AN_n^{0},AN_n^{1},\ldots,AN_n^{n-1}$.
For each $i\in[n]$, $AN_n^{i}$ is isomorphic to $AN_{n-1}$.
By Theorem 1 in \cite{Z1}, $AN_{n}$ is $(n-1)$-regular and $(n-1)$-connected.

\begin{lem}{\rm(\cite{H})}\label{AN-1}
Let $AN_n$ be the alternating group network for $n\geq 3$ .
\begin{enumerate}
\item [{\rm (1)}] Each vertex in $AN_{n}$ has exactly one extra neighbor.
\item [{\rm (2)}] $AN_{n}$ has no $4$-cycle and $5$-cycle.
\item [{\rm (3)}] Let $u$ and  $v$ be any two distinct vertices of $AN_n$, then $cn(AN_{n}:u,v)\leq 1$.
\end{enumerate}
\end{lem}

\begin{lem}{\rm(\cite{Z99})}\label{ANn4}
Let $F$ be a vertex-cut of $AN_n$ for $n\geq 5$. If $|F|\leq 2n-5$, then $AN_n-F$ satisfies one of the following conditions:
\begin{enumerate}
\item [{\rm (1)}] $AN_n-F$ has two components, one of which is a trivial component.
\item [{\rm (2)}] $AN_n-F$ has two components, one of which is an edge.
Moreover, if $|F|=2n-5$, $F$ is formed by the neighbor of the edge.
\end{enumerate}
\end{lem}

\subsection{BC networks}

\begin{defi}
The $1$-dimensional BC network $X_{1}$ is a complete
graph with two vertices. The {\it $n$-dimensional BC
network} $X_{n}$ is defined as follows: $V(X_{n})=V(G_{1})\cup V(G_{2})$
and $E(X_{n})=E(G_{1})\cup E(G_{2})\cup M$, where $G_{1}, G_{2}\in L_{n-1}$,
and $M$ is a perfect matching between $V(G_{1})$ and $V(G_{2})$, where $L_{k}=\{X_{k}: X_{k} \ {\rm is\ an}\ k$-${\rm dimensional\ BC\ network}\}$.
\end{defi}

\begin{lem}\label{bc2}
Let $G=X_n\in L_n$ for $n\geq 1$. Then
\begin{enumerate}
\item [{\rm (1)}] {\rm(\cite{FH})} $|V(G)|=2^n$, $|E(G)|=n\cdot2^{n-1}$, $G$ is $n$-regular and triangle-free.
\item [{\rm (2)}] {\rm(\cite{FH},\cite{VR})} $\kappa(G)=n$.
\item [{\rm (3)}] {\rm(\cite{Z08})} $cn(G)=2$.
\end{enumerate}
\end{lem}

\begin{lem}{\rm(\cite{Z08})}\label{bc4}
For any $X_n\in L_n$, let $F\subseteq V(X_n)$ with $|F|\leq 2n-3$ be a vertex-cut of $X_n$.
Then $X_n-F$ has two components, one of which is a trivial component.
\end{lem}

\subsection{The $k$-ary $n$-cube networks}

\begin{defi}
The {\it $k$-ary $n$-cube}, denoted by $Q_{n}^{k}$, where $k\geq2$ and $n\geq1$ are integers, is a graph consisting of $k^{n}$ vertices, each of these vertices has the form
$u=u_{n-1}u_{n-2}\cdots u_{0}$, where $u_{i}\in \{0,1,\ldots,k-1\}$ for $0\leq i\leq n-1$. Two vertices $u=u_{n-1}u_{n-2}\cdots u_{0}$ and $v=v_{n-1}v_{n-2}\cdots v_{0}$ in $Q_{n}^{k}$
are adjacent if and only if there exists an integer $j$, where $0\leq j\leq n-1$, such that $u_{j}=v_{j}\pm 1($mod $k)$ and $u_{i}=v_{i}$ for every $i\in \{0,1,\ldots,n-1\}\setminus\{j\}$.
In this case, $(u, v)$ is a {\it j-dimensional edge}.
\end{defi}

For convenience, $``($mod $k)"$ does not appear in similar expressions in the remainder of the paper.
Note that each vertex has degree $2n$ for $k\geq 3$ and has degree $n$ for $k=2$. Clearly, $Q_{1}^{k}$ is a cycle of length $k$,
$Q_{n}^{2}$ is an {\it $n$-dimensional hypercube}, $Q_{2}^{k}$ is a $k\times k$ {\it wrap-around mesh}.

$Q_{n}^{k}$ can be partitioned over the $j$th-dimension, for a $j\in [n-1]$, into $k$ disjoint subcubes, denoted by $Q_{n-1}^{k}[0],Q_{n-1}^{k}[1],\ldots,Q_{n-1}^{k}[k-1]$,
by deleting all the $j$-dimensional edges from $Q_{n}^{k}$. For convenience, abbreviate these as $Q[0],Q[1],\ldots, Q[k-1]$ if there is no ambiguity.
Moreover, $Q[i]$ for $0\leq i\leq k-1$ is isomorphic to the $k$-ary $(n-1)$-cube. For each vertex $u \in V(Q[i])$, the neighbor which is not in $V(Q[i])$ is called the {\it extra neighbor}.
For $i\in [k-1]$, $u \in V(Q[i])$, the two extra neighbors of $u$ are in different subgraphs $Q[i+1]$ and $Q[i-1]$, respectively.

\begin{lem}\label{Qnk2}
Let $Q_{n}^{k}$ be a $k$-ary $n$-cube, where $k\geq2$ and $n\geq 1$ are integers.
\begin{enumerate}

\item [{\rm (1)}] {\rm(\cite{D97})}  $Q_{n}^{k}$ is $2n$-regular and $2n$-connected for $k\geq 3$
 and $n$-regular and $n$-connected for $k=2$.
\item [{\rm (2)}] {\rm(\cite{D,G,Hsieh})} For any $x,y\in V(Q_{n}^{k})$, $k\geqslant 2$,
 $$
 |N(x)\cap N(y)|=\left\{\begin{array}{ll}
 1& \ \text{if $xy\in E(Q_{n}^k)$ and $k=3$};\\
 2& \ \text{if $xy\notin E(Q_{n}^k)$ and $N(x)\cap N(y)\ne\emptyset$};\\
 0& \ \text{otherwise.}
 \end{array}\right.
 $$
 \end{enumerate}
\end{lem}

\begin{lem}\label{Qnk4}
\begin{enumerate}
\item [{\rm (1)}]{\rm(\cite{E89})} If $F\subseteq V(Q_{n}^{2})$ with $|F|\leq 2n-3$ is a vertex cut of $ Q_{n}^{2}$ for $n\geq 2$, then $Q_{n}^{2}-F$ has two components, one of which is a trivial component.
\item [{\rm (2)}] {\rm(\cite{D,G})}
If $F\subseteq V(Q_{n}^{3})$ with $|F|\leq 4n-4$ is a vertex cut of $ Q_{n}^{3}$ for $n\geq 2$, then $Q_{n}^{3}-F$ has two components, one of which is a trivial component.
\item [{\rm (3)}] {\rm(\cite{D,G15})}
If $F\subseteq V(Q_{n}^{k})$ is a vertex cut of $Q_{n}^{k}$ with $|F|\leq 4n-3$ for $n\geq 2$ and $k\geq4$, then $Q_{n}^{k}-F$ has two components, one of which is a trivial component.

\end{enumerate}
\end{lem}

\subsection{Split-star networks $S_{n}^{2}$}

Cheng et al.~\cite{CLP} propose the Split-star networks as alternatives to the star graphs and companion graphs with the alternating group graphs.

\begin{defi}
Given two positive integers $n$ and
$k$ with $n>k$, note that $[n]=\{1,2,\ldots, n\}$, and let $\mathcal{P}_n$ be a set of $n!$ permutations on $[n]$. The {\it $n$-dimensional Split-star
network}, denoted by $S_n^2$, such that $V(S_n^2)=\mathcal{P}_n$, $E(S_n^2)=\{(p,q)| \ p$ (resp. $q$) can be obtained from $q$
(resp. $p$) by either a $2$-exchange or a $3$-rotation $\}$. Where
\begin{enumerate}
\item [{\rm (1)}] A {\it $2$-exchange interchanges} the symbols in 1st position
and 2nd position.
\item [{\rm (2)}] A {\it $3$-rotation rotates} the symbols in three positions
labeled by the vertices of a triangle %in Fig. 1,
in which three vertices of the triangle are $1, 2$ and $k$ for some
$k\in \{3, 4,\ldots, n\}$.
\end{enumerate}
\end{defi}

Let $V_n^{n:i}$ be the set of all vertices in $S_n^2$ with the $n$th position having value $i$, i.e.,
$V_n^{n:i}=\{p|p=x_1x_2\cdots x_{n-1}i$, $x_j\in \{1,2,\ldots,i-1,i+1,\ldots n\}$
($1\leq j\leq n-1$) are do not care symbols $\}$.
The set $\{V_n^{n:i}|1\leq i\leq n\}$ forms a partition $V(S_n^2)$.
Let $S_n^{2:i}$ denote the subgraph of $S_n^2$ induced by $V_n^{n:i}$, i.e., $S_n^{2:i}=S_n^2[V_n^{n:i}]$.
It is easy to know that $S_n^{2:i}$ is isomorphic to $S_{n-1}^2$.
Every vertex $v\in S_n^{2:i}$ has exactly two neighbors, called extra neighbors, outside of $S_n^{2:i}$;
 moreover these two neighbors belong to different $S_n^{2:j}$s where $j\neq i$.
We call these neighbors as the extra neighbors of $v$. We call these edges, whose end-vertices belong to different subgraphs, as {\it cross edges}. Let $S_{n,E}^2$ be a subgraph of $S_n^2$ induced by the set of even permutations, in which the adjacency rule is precisely the $3$-rotation. We know that $S_{n,E}^2$ is the alternating group graph $AG_n$~\cite{jhd93}.
Let $S_{n,O}^2$ be a subgraph of $S_n^2$ induced by the set of odd permutations, in which the adjacency rule is precisely the $3$-rotation. We have that $S_{n,O}^2$ is also isomorphic to $AG_n$ and $S_{n,O}^2$ is isomorphic $S_{n,E}^2$ via the $2$-exchange
$\phi(a_1a_2a_3\cdots a_n)=a_2a_1a_3\cdots a_n$. Hence, there are $\frac{n!}{2}$ matching edges between $S_{n,O}^2$ and $S_{n,E}^2$. Indeed, the Split-star network $S_n^2$ is introduced in~\cite{CLP01} which is the companion graph of $AG_n$.

\begin{lem}{\rm(~\cite{CL,CLP,CLP01})}\label{Sn2}
Let $S_{n}^{2}$ be the $n$-dimensional split-star network.
\begin{enumerate}
\item [{\rm (1)}] $S_{n}^{2}$ is $(2n-3)$-regular and $\kappa(S_{n}^{2})=2n-3$ for $n\geq 2$.
\item [{\rm (2)}] Two extra neighbors of every vertex in $S_n^{2:i}$ are in distinct induced subgraphs and these two extra neighbors are adjacent. For any two vertices in the same subgraph $S_n^{2:i}$, their extra neighbors in other subgraphs are different. There is one to one correspondence between the subgraph $S_{n,O}^2$ and the subgraph $S_{n,E}^2$.
\item [{\rm (3)}] Let $x,y$ be any two vertices of $S_{n}^{2}$, then
$$
 |N(x)\cap N(y)|=\left\{\begin{array}{ll}
 1& \ \text{if $d(x,y)=1$};\\
 2& \ \text{if $d(x,y)=2$};\\
 0& \ \text{if $d(x,y)\geq 3$.}
 \end{array}\right.
 $$
\end{enumerate}
\end{lem}

\begin{lem}{\rm(~\cite{LX15})}\label{Sn24}
If $F\subseteq V(S_n^2)$ with $|F|\leq 4n-10$ is a vertex cut of $S_n^2$ for $n\geq 4$, then $S_n^2-F$ has two components, one of which is a trivial component.
\end{lem}

\subsection{Cayley graphs generated by transposition trees $\Gamma_n$}

Note that $\mathcal{P}_n$ is a group of all permutations on $[n]$. For convenience,
$(ij)$, which is called a {\it transposition}, denotes the permutation that swaps the elements at position $i$ and $j$, that is $(ij)p_1p_2\ldots p_i\ldots p_j\ldots p_n=p_1p_2\ldots p_j\ldots p_i\ldots p_n$.

\begin{defi}
Let $\mathcal{P}_n$ be symmetric group on $[n]$, and the generating
set $S$ to be a set of transpositions. A graph $G(S)$ with vertex
set $[n]$, where there is an edge between $i$ and $j$ if and only if the transposition $(ij)$ belongs to $S$, is called the {\it transposition generating graph}.
When $G(S)$ is a tree, we call $G(S)$ a {\it transposition tree}.
The Cayley graphs $Cay(\mathcal{P}_n, S)$ obtained by transposition trees are called {\it Cayley graphs generated by transposition trees}, denoted by $\Gamma_n$.
\end{defi}

If $G(S)\cong K_{1,n-1}$, $Cay(\mathcal{P}_n, S)$ is called the \emph{star graph}, denoted by $S_n$. If $G(S)\cong P_n$, that is the transposition tree is a path $P_n$ with $n$ vertices, then $Cay(\mathcal{P}_n, S)$ is called the \emph{bubble-sort graph}, denoted by $B_n$.

Let $\Gamma_n^i$ be the subgraph of $\Gamma_n$ spanned by vertices corresponding to permutations with $i$ in the last position. Then $\Gamma_n$ can be divided into $n$ subgraphs $\Gamma_{n-1}^1$, $\Gamma_{n-1}^2$, $\cdots$, $\Gamma_{n-1}^n$ and each $\Gamma_{n-1}^i$ is isomorphic to $\Gamma_{n-1}$ for $i\in[n]$. For $u\in V(\Gamma_{n-1}^{i})$, denoted by $u'=u(1n)$ the unique neighbor of $u$ outside $\Gamma_{n-1}^{i}$, called the extra neighbor of $u$.

\begin{lem}\label{Gamma}
Let $\Gamma_n$ be the Cayley graphs generated by transposition trees for $n\geq 3$.
\begin{enumerate}
\itemsep -1pt
\item [{\rm (1)}] {\rm(\cite{C08})} $\kappa(\Gamma_n)=n-1$.
\item [{\rm (2)}] {\rm(\cite{C08})} $\Gamma_n$ has the girth $4$ unless $\Gamma_n$ is the star graph which has girth $6$. $\Gamma_n$ does not have $K_{2,3}$ as a subgraph.
\item [{\rm (3)}] {\rm(\cite{Y10})} For any two distinct vertices $u,v\in \Gamma_n$, $|N_{\Gamma_n}(u)\cap N_{\Gamma_n}(v)|=1$ if $\Gamma_n=S_n$; Otherwise $|N_{\Gamma_n}(u)\cap N_{\Gamma_n}(v)|\leq 2$.
\end{enumerate}
\end{lem}

\begin{lem}{\rm(\cite{C08,Y10})}\label{Gamma4}
If $F\subseteq V(\Gamma_n)$ with $|F|\leq 2n-5$ is a vertex cut of $\Gamma_n$ for $n\geq 4$, then $\Gamma_n-F$ has two components, one of which is a trivial component.
\end{lem}

\subsection{Cayley graphs generated by $2$-trees}

\begin{defi}
Let $\Gamma$ be the alternating
group, the set of even permutations on $\{1, 2,\ldots,n\}$,
and the generating set $\Delta$ to be a set of $3$-cycles.
To get an undirected Cayley graph, we will assume
that whenever a $3$-cycle $(abc)$ is in $\Delta$, so is its inverse,
$(acb)$. Since $(abc)$, $(bca)$ and $(cab)$ represent the
same permutation, the set $\{a,b,c\}$ uniquely represents this
$3$-cycle and its inverse. So we can depict $\Delta$ via a hypergraph
with vertex set $[n]$, where a hyperedge of size $3$ corresponds
to each pair of a $3$-cycle and its inverse in $\Delta$.
\end{defi}

It is easy to see that the Cayley graph generated by the
$3$-cycles in $\Delta$ is connected if its corresponding hypergraph $H$ is
connected. Since an interconnection network needs to be
connected, we require $H$ graph to be connected.

In general,
this graph may have extra $K_{3}$'s formed by vertices
that do not correspond to a $3$-cycle in $\Delta$. We will avoid
this possibility by considering a simpler case when $H$ has
a tree-like structure. Such a graph is built by the following
procedure. We start from $K_{3}$, then repeatedly add a new
vertex, joining it to exactly two adjacent vertices of the previous
graph. Any graph obtained by this procedure is called
a $2$-tree. If $v$ is a vertex of a $2$-tree $H$ with the property
that $H$ can be generated in such a way that $v$ is the last
vertex added, then $v$ is called a leaf of the $2$-tree.

The {\it alternating group graph} $AG_{n}$ \cite{J}, can be viewed as the
Cayley graph generated by the graph  having a
tree-like (in fact, star-like) structure of triangles.

It is easy to prove that if two $2$-trees are isomorphic, then the corresponding Cayley graphs will also be isomorphic;
hence without loss of generality we may assume that vertex n is the tail of the $2$-tree.
For $n\geq4$, the vertices corresponding to even permutations ending with $i$
induce a subgraph $\Gamma_{n-1}^{i}(\Delta)$ that is also a Cayley graph generated by a
2-tree $\Delta'$, which is obtained by deleting the edges corresponding to the two $3$-cycles in $\Delta$ containing $n$.
Thus we obtain the following result of the recursive structure of $\Gamma_{n}(\Delta)$:

\begin{lem}{\rm(\cite{C1})}\label{prop-1}
Let $\Gamma_{n}(\Delta)$ be a Cayley graph generated by the $2$-tree $\Delta$, $\Delta'=\Delta-\{n\}$, $n\geq 4$.
Then
 \begin{enumerate}
  \item [{\rm (1)}] $\Gamma_{n}(\Delta)$ consists of $n$ vertex-disjoint subgraphs, $\Gamma_{n-1}^{1}(\Delta),\Gamma_{n-1}^{2}(\Delta),\ldots,\Gamma_{n-1}^{n}(\Delta)$, each isomorphic to $\Gamma_{n-1}(\Delta')$.
  \item [{\rm (2)}] $\Gamma_{n-1}^{i}(\Delta)$ has $(n-1)!/2$ vertices, and it is $(2n-6)$-regular for all $i$.
  \item [{\rm (3)}] There are exactly $(n-2)!$ independent edges between $\Gamma_{n-1}^{i}(\Delta)$ and $\Gamma_{n-1}^{j}(\Delta)$ for all $i\neq j$.
  \item [{\rm (4)}] Each vertex in $\Gamma_{n-1}^{i}(\Delta)$ has exactly two neighbors outside $\Gamma_{n-1}^{i}(\Delta)$; these two outside neighbors are in different $\Gamma_{n-1}^{k}(\Delta)$'s, and there is an edge between them. Thus every vertex forms a triangle with its two outside neighbors.
  \item [{\rm (5)}] $\Gamma_{n}(\Delta)$ does not contain $K_{4}-e$, that is, $K_{4}$
with an edge deleted, and $K_{2,3}$ as a subgraph. For any two vertices $u$ and $v$, $|N(u)\cap N(v)|=1$ if $d(u,v)=1$, $|N(u)\cap N(v)|\leq 2$ otherwise.
  \end{enumerate}

\end{lem}

\begin{lem}{\rm(\cite{C})}\label{2-con}
Let $G=\Gamma_{n}(\Delta)$ be a Cayley graph generated by the $2$-tree $\Delta$ for $n\geq 4$. Then $G$ is maximally
connected, i.e., $G$ is $(2n-4)$-regular and $(2n-4)$-connected.
\end{lem}

\begin{lem}{\rm(\cite{C})}\label{Delta}
Let $G=\Gamma_{n}(\Delta)$ be a Cayley graph generated by the $2$-tree $\Delta$ for $n\geq 4$,
and let $T$ be a set of vertices in $G$ such that $|T|\leq 4n-11$. If $n\geq 5$, then $G-T$ satisfies one of the following conditions:
 \begin{enumerate}
  \item [{\rm (1)}] $G-T$ is connected.
  \item [{\rm (2)}] $G-T$ has two components, one of which is a singleton.
  \item [{\rm (3)}] $G-T$ has two components, one of which is a $K_{2}$. Moreover, $|T|=4n-11$, and the set $T$ is formed by the neighbors of the two vertices in the $K_{2}$.
 \end{enumerate}
When $n=4$, there are two additional possibilities. In both cases, $G-T$ has two components, one of which is a $4$-cycle. The other component is either a $4$-cycle if $|T|=4$ or a path with $3$ vertices if $|T|=5$.
\end{lem}

\subsection{Burnt pancake networks $BP_n$}

Gates and Papadimitriou~\cite{GP79} introduced the burnt pancake problem in 1979. Burnt pancake problem relates to the construction of networks of parallel processors.

Let $n$ be a positive integer. We use $[n]$ to denote the set $\{1,2,\ldots,n\}$. To save space, the negative sign may be placed on the top of an expression. Thus, $\bar{i}=-i$. We use $\langle n\rangle$ to denote the set $[n] \cup \{\bar{i}|i\in [n]\}$. A {\it signed permutation} of
$[n]$ is an $n$-permutation $u_1u_2\cdots u_n$ of $\langle n\rangle$ such that
$|u_1||u_2|\cdots |u_n|$ taking the absolute value of each element, forms a permutation of $[n]$.
For a signed permutation $u=x_1x_2\cdots x_i\cdots x_n$ of $\langle n \rangle$, the $i$-th prefix reversal of u,
denoted by $u^i$ is $u^i=\bar{x}_i\bar{x}_{i-1}\cdots \bar{x}_1x_{i+1}\cdots x_n, 1\leq i\leq n$.
For example, let $u=1\bar{2}4\bar{3}5$; then $u$ is a signed permutation of $[5]$,
$u^2=2\bar{1}4\bar{3}5$, $u^5=\bar{5}3\bar{4}2\bar{1}$.
%(2) = 21435 and u
%(5) = 53421.

\begin{defi}
An {\it $n$-dimensional burnt pancake network} $BP_n$ is
defined to be an $n$-regular graph $G$ with $n!2^n$ vertices, each of which has a unique label from
the signed permutation of $\langle n \rangle$.
Two vertices $u$ and $v$ are adjacent in $BP_n$ if and only if $u^i=v$ for some unique $i$ ($1\leq i\leq n$).
Such an edge $uv$ is called an {\it $i$-dimensional edge} and $v$ is called the {\it $i$-neighbor} of $u$.
It is seen that every vertex has a unique $i$-neighbor for $1\leq i\leq n$.

\end{defi}

\begin{lem}{\rm(\cite{CWHC09,C2011,IK10})}\label{BPn}
An $n$-dimensional burnt pancake network $BP_n$ has the following combinatorial properties.
\begin{enumerate}
  \item [{\rm (1)}] $BP_n$ is $n$-regular with $n!\times 2^n$ vertices and $n!\times 2^{n-1}$ edges.
  \item [{\rm (2)}] $\kappa(BP_n)=n$, the girth of $BP_n$($n\geq 3$) is $g(BP_n)=8$.
  \item [{\rm (3)}] $BP_n$ can be decomposed into $2n$ vertex-disjoint subgraphs, denoted $BP_n^i$, by fixing the
  symbol in the last position $n$, in which the symbol in the $n$th position is $i$, where $i\in [n]$. Obviously,
  $BP_n^i$ is isomorphic to $BP_{n-1}$. The number of cross edges between any two subgraphs, $BP_n^i$ and $BP_n^j$
  ($i\neq j, i,j\in[n]$), is $|E(i,j)|=(n-2)!\times2^{n-2}$ if $i\neq \bar{j}$; otherwise, $|E(i,j)|$=0.
  For a vertex $v\in V(BP_n^i)$, $v$ has exactly one neighbor outside $BP_n^i$, called the {\it extra neighbor} of $v$.
\end{enumerate}
\end{lem}

\begin{lem}{\rm(\cite{SXZC15})}\label{BPn4}
For any subset $F\subseteq V(BP_n)$ with $|F|\leq 2n-2$ is a vertex-cut of $BP_n$ for $n\geq 4$, then $BP_n-F$ satisfies
one of the following conditions.
\begin{enumerate}
  \item [{\rm (1)}] $BP_n-F$ has two connected components, one of which is a trivial component;
  \item [{\rm (2)}] $BP_n-F$ has two connected components, one of which is an edge. Furthermore,
  $F$ is the neighborhood of this edge with $|F|=2n-2$.
\end{enumerate}
\end{lem}

\section{Main result}

\f In this section, the relationship between the pessimistic diagnosability
under the PMC model and the
extra connectivity with some restricted conditions will
be proposed.

\begin{lem}\label{lem-k}
Let $G$ be a $k$-regular graph.
Let $u$ and $v$ be two distinct vertices
in $G$, if  $cn(G;u,v)\leq 2$, then $|N_{G}(\{u,v\})|\geq 2k-2-l$, where $l=l(G)=\max \{cn(G;u,v): (u,v)\in E(G)\}$, i.e., $l=l(G)$ be the maximum number of common neighbors between any two adjacent vertices.
\end{lem}

\f\demo Since $cn(G;u,v)\leq 2$, if $u$ is non-adjacent to $v$,
then $|N_{G}(\{u,v\})|=|N_{G}(u)|+|N_{G}(v)|- cn(G;u,v)\geq 2k-2\geq 2k-2-l$.
Otherwise, $u$ is adjacent to $v$, $|N_{G}(\{u,v\})|=|N_{G}(u)|-1+|N_{G}(v)|-1-cn(G;u,v)\geq 2(k-1)-l$.
As a result, $|N_{G}(\{u,v\})|\geq 2k-2-l$.
\hfill\qed

Tsai and Chen~\cite{T} derived the following result which characterizes a graph for $t/t$-diagnosability.

\begin{lem}{\rm(\cite{T})}\label{lem-5}
A graph $G$ is $t/t$-diagnosable if and only if for each vertex set $S\subseteq V(G)$ with $|S|=p$, $0\leq p\leq t-1$,
$G-S$ has at most one trivial component and each nontrivial component $C$ of $G-S$ satisfies $|V(C)|\geq 2(t-p)+1$.
\end{lem}

The following result is useful.

\begin{lem}{\rm(\cite{F})}\label{lem-3}
Let $G$ be a connected graph and $U\subseteq V(G)$. Then, $|N_{V(G)-U}(U)|\geq \kappa(G)$ if $|V(G)-U|\geq \kappa(G)$, otherwise, $|N_{V(G)-U}(U)|=|V(G)-U|$.
\end{lem}

\begin{theorem}\label{main}
Let $G$ be a $k$-regular $k$-connected ($k\geq 5$) graph with order $N$.
Let $U$ be a subset of $V(G)$ and $l=l(G)$ be the maximum number of common neighbors between any two adjacent vertices.
Suppose further that all of the following conditions hold:
\begin{enumerate}
\item [{\rm (1)}] $N\geq 4k-2$.
\item [{\rm (2)}] $cn(G)\leq 2$.
\item [{\rm (3)}] If $2\leq |U|\leq 2(2k-4-l)$, then $|N_{G}(U)|\geq 2k-2-l$.
\item [{\rm (4)}] Let $F\subseteq V(G)$ be a vertex-cut of $G$.
If $|F|\leq 2k-3-l$, then $G-F$ has a large component and a small component which is a trivial component.
\end{enumerate}
Then, $t_{p}(G)=2k-2-l=\kappa_1(G)$.
\end{theorem}

\f\demo We first prove $t_{p}(G)\leq 2k-2-l$. Suppose $t_{p}(G)\geq 2k-2-l+1$,
then $G$ is $(2k-2-l+1)/(2k-2-l+1)$-diagnosable.
Let $(u,v)$ be an edge of $G$ such that $|N_{G}(u)\cap N_{G}(v)|=l$.
Let $S=N_{G}(\{u,v\})$. Then $|S|=2k-2-l\leq t_{p}(G)-1$. An edge $\{u,v\}$ is a connected component of $G-S$, say $C$.
By Lemma~\ref{lem-5}, $|V(C)|\geq 2(t_{p}(G)-|S|)+1\geq 2[(2k-2-l+1)-(2k-2-l)]+1=3$, which is a contradiction. Thus, $t_{p}(G)\leq 2k-2-l$.

Secondly, we show $t_{p}(G)\geq2k-2-l$, i.e., $G$ is $(2k-2-l)/(2k-2-l)$-diagnosable. Suppose $G$ is not $(2k-2-l)/(2k-2-l)$-diagnosable,
by Lemma~\ref{lem-5}, there exists a vertex set
$S\subseteq V(G)$ with $|S|=p$, $0\leq p\leq 2k-3-l$ such that $G-S$ contains more than one trivial components or
contains a nontrivial component $C$ with $|V(C)|\leq 2(2k-2-l-p)$. The following cases should be considered.

Case 1. $G-S$ contains more than one trivial components.

Suppose $C_{1}=\{u\}$ and $C_{2}=\{v\}$ are two distinct trivial components of $G-S$. By Condition (2) and
Lemma~\ref{lem-k}, $|N_{G}(\{u,v\})|\geq 2k-2-l$.
Note that $N_{G}(\{u,v\})\subseteq S$, this implies that $|S|\geq 2k-2-l$, which is a contradiction.

Case 2. $G-S$ contains a nontrivial component $C$ with $2\leq |V(C)|\leq 2(2k-2-l-p)$.

Suppose $p\leq 1$. Since the connectivity of $G$ is $k\geq 5>p$, $G-S$ is connected.
It implies $C=G-S$. By $|V(C)|=|V(G)|-|S|=N-p\geq N-1$, Condition $(1)$ and $l\leq cn(G)\leq 2$,
one has $|V(C)|\geq 4k-3\geq 2(2k-2-l-p)+1$ which is a contradiction.

Now consider $2\leq p\leq 2k-3-l$. Since $2\leq |V(C)|\leq 2(2k-2-l-p)$, so $2\leq |V(C)|\leq 2(2k-4-l)$.
By condition $(3)$, $|N_{G}(V(C))|\geq 2k-2-l$. Since $C$ is a connected component of $G-S$, $N_{G}(V(C))\subseteq S$. This implies $p=|S|\geq 2k-2-l$,
which is a contradiction for the fact that $p=|S|\leq 2k-3-l$. Thus, $t_{p}(G)\leq 2k-2-l$.

Next we prove $2k-2-l=\kappa_1(G)$. Let $(u,v)$ be an edge of $G$ such that $|N_{G}(u)\cap N_{G}(v)|=l$.
Let $S=N_{G}(\{u,v\})$. Then $|S|=2k-2-l$. If $G-S=\{(u,v)\}$, then $|V(G)|=|S|+2=2k-l<4k-2$
for $k\geq 5$ which contradicts with Condition $(1)$. If $G-S$ has a trivial component which contains only one vertex, say $\{x\}$, then $G-S$ has at least
two components: $\{x\}$ and the edge $(u,v)$. By $cn(G)\leq 2$,
then $|S|\geq 2k-2-l+(k-4)=3k-6-l$. Note $3k-6-l>2k-2-l$ for $k\geq 5$, it is a contradiction.
Thus, $G-S$ has no trivial component, i.e., $S$ is an extra vertex cut of $G$, which implies $\kappa_1(G)\leq 2k-2-l$.
On the other hand, by condition $(4)$, $\kappa_1(G)\geq 2k-2-l$. Thus, $\kappa_1(G)=2k-2-l$.

By above discussion, $t_{p}(G)=2k-2-l=\kappa_1(G)$. \hfill\qed

\section{Application to some interconnection networks}

As applications of Theorem \ref{main}, in this section, we
determine the pessimistic diagnosability and extra connectivity
for some well-known interconnection networks,
including the alternating group graph $AG_n$, the alternating group network $AN_n$,
the $k$-ary $n$-cube networks $Q_n^k$, BC networks $X_n$, split-star networks $S_{n}^{2}$,
Cayley graphs generated by transposition trees $\Gamma_n$,
Cayley graphs generated by $2$-trees, burnt pancake networks $BP_n$.

\subsection{Application to the alternating group graphs $AG_n$}

\begin{rmk}
It is known that $\kappa_{1}(AG_n)=4n-11$ for $n\geq 5$
determined by Lin {\it et al.}~\cite{LZX}
and $t_p(AG_n)=4n-11$ obtained by Tsai~\cite{T1}.
As a corollary of Theorem~\ref{main}, we immediately obtain the following result which contains the above result.
\end{rmk}

\begin{cor}\label{cor1}
Let $AG_n$ be the $n$-dimensional alternating group graph for $n\geq 5$. Then $t_p(AG_n)=4n-11=\kappa_{1}(AG_n)$.
\end{cor}

\f\demo Obviously, $N=|V(AG_n)|=\frac{n!}{2}$, $k=2n-4\geq 6$ for $n\geq 5$, $l=l(AG_n)=1$.

Note that $N=\frac{n!}{2}\geq 4(2n-4)-2$ for $n\geq 5$,
Conditions $(1)$ in Theorem~\ref{main} holds.
Conditions $(2)-(4)$ in Theorem~\ref{main} hold
by Lemmas~\ref{AGn1},~\ref{AGn3} and~\ref{AGn4}, respectively.
Thus, $AG_n$ satisfies all conditions in Theorem~\ref{main}, $t_p(AG_n)=4n-11=\kappa_{1}(AG_n)$ for $n\geq 5$. \hfill\qed

\subsection{Application to the alternating group networks}

Zhou~\cite{Z99} derived $\kappa_1(AN_n)=2n-5$ for $n\geq 4$.
However, $t_p(AN_n)$ has not been determined so far.
We can deduce the result as a corollary of Theorem~\ref{main} as following.
Notice that for $AN_n$, $k=n-1$, $l=1$ in Theorem~\ref{main}.

\begin{lem}\label{ANn3}
Let $AN_n$ be the $n$-dimensional alternating group network for $n\geq 4$.
If $U$ is a subset of $V(AN_{n})$ and $2\leq |U|\leq 2(2k-4-l)=4n-14$, then $|N_{AN_{n}}(U)|\geq 2n-5$.
\end{lem}

\f\demo The Lemma can be proved by using the induction on $n$.
It is easy to verify that $|N_{AN_{4}}(U)|\geq 3$ for $|U|=2$ by Lemma~\ref{lem-k}.
We assume that the lemma is true for $AN_{m}$, where $m$ is an integer with $5\leq m\leq n-1$, we will prove the result for $AN_{n}$.

Recall that $AN_{n}$ is constructed by $n$ disjoint $AN_{n-1}$'s, denoted by $AN_n^{i}$ for $i\in [n]$.
Let $U_{i}=U\cap V(AN_n^{i})$ and $\overline{AN_{n}^{i}}=AN_{n}-AN_{n}^{i}$ for $i\in [n]$. Without loss of generality, we may assume that $|U_{1}|\geq |U_{2}|\geq \ldots \geq |U_{n}|$.
The following cases should be considered.

Case 1. $|U_{1}|\leq 1$.

In this case, $|U_{i}|\leq 1$ for all $i\in [n]$. Clearly, $2\leq |U|\leq n$ because of $i\leq n$.
The Lemma follows if $|U|=2$ by Lemma~\ref{lem-k}.
Now assume that $3\leq |U|\leq n$. Since $AN_{n}$ is $(n-1)$-regular and $AN_n^{i}$ is isomorphic to $AN_{n-1}$, $|N_{AN_n}(U)|\geq 3\kappa(AN_n^{i})=3(n-2)\geq 2n-5$ for $n\geq 7$.

Case 2. $2\leq |U_{1}|\leq 4n-19$.

By inductive hypothesis in $AN_n^{1}$, $|N_{AN_n^{1}}(U_{1})|\geq 2(n-1)-5=2n-7$.
%Since $|U|\leq 4n-14$ and $|U_{1}|\geq |U_{2}|\geq \ldots \geq |U_{n}|$, $|U_{2}|\leq 2n-7$.
If $U=U_{1}$, $|N_{AN_n}(U)|=|N_{AN_n^{1}}(U_{1})|+|N_{\overline{AN_{n}^{1}}}(U_{1})|\geq 2n-7+|U_1|\geq 2n-5$. Assume $U\neq U_1$ in the following.
If $|U_{2}|=1$, $|N_{AN_n^{2}}(U_{2})|=\kappa(AN_n^{2})=n-2$.
Note that $AN_n^{1}$ and $AN_n^{2}$ are vertex disjoint, $|N_{AN_n}(U)|\geq |N_{AN_n^{1}}(U_{1})|+|N_{AN_n^{2}}(U_{2})|\geq 3n-9\geq 2n-5$ for $n\geq 5$.
Now consider $2\leq |U_{2}|\leq |U_1|\leq 4n-19$, by inductive hypothesis in $AN_n^{2}$, $|N_{AN_n^{2}}(U_{2})|\geq 2(n-1)-5=2n-7$.
Thus, $|N_{AN_n}(U)|\geq |N_{AN_n^{1}}(U_{1})|+|N_{AN_n^{2}}(U_{2})|\geq4n-14 \geq 2n-5$ for $n\geq 5$.

Case 3. $4n-18\leq|U_{1}|\leq 4n-14$.

Since the connectivity of $AN_{n}^{1}$ is $n-2$, and $\frac{(n-1)!}{2}-|U_{1}|\geq n-2=\kappa(AN_{n}^{1})$ for $n\geq 5$,  by Lemma~\ref{lem-3}, $|N_{AN_n^{1}}(U_{1})|\geq n-2$.
By Lemma~\ref{AN-1}, $|N_{\overline{AN_{n}^{1}}}(U_{1})|=|U_{1}|$.
If $U=U_{1}$, $|N_{AN_n}(U)|\geq |N_{AN_n^{1}}(U_{1})|+|N_{\overline{AN_{n}^{1}}}(U_{1})|\geq (n-2)+4n-18=5n-20\geq 2n-5$ for $n\geq 5$.
In the following, we assume the case of $U\neq U_{1}$.
Note that $U\neq U_{1}$ and $|U-U_{1}|\leq 3$, so $1\leq |U_{2}|\leq 3$.

If $|U_{2}|=1$, recall that $AN_{n}$ is $(n-1)$-regular and $AN_n^{i}$ is isomorphic to $AN_{n-1}$, $|N_{AN_n^{2}}(U_{2})|=\kappa(AN_n^{2})=n-2$.
Hence, $|N_{AN_n}(U)|\geq |N_{AN_n^{1}}(U_{1})|+|N_{AN_n^{2}}(U_{2})|\geq 2n-4\geq 2n-5$ for $n\geq 5$.
Now suppose that $2\leq |U_{2}|\leq 3$. Since $\frac{(n-1)!}{2}-|U_{2}|\geq n-2=\kappa(AN_{n}^{2})$ for $n\geq 5$,
by Lemma~\ref{lem-3}, $|N_{AN_n^{2}}(U_{2})|\geq n-2$.
Thus, $|N_{AN_n}(U)|\geq |N_{AN_n^{1}}(U_{1})|+|N_{AN_n^{2}}(U_{2})|\geq 2(n-2)\geq 2n-5$ for $n\geq 5$.
%Since $3\leq 4n-19$ for $n\geq 6$, by inductive hypothesis in $AN_n^{2}$, $|N_{AN_n^{2}}(U_{2})|\geq 2(n-1)-5=2n-7$.
%Thus, $|N_{AN_n}(U)|\geq |N_{AN_n^{1}}(U_{1})|+|N_{AN_n^{2}}(U_{2})|\geq (n-2)+(2n-7)\geq 2n-5$ for $n\geq 6$.
%If $n=5$,

By the above cases, the Lemma holds.\hfill\qed

\begin{cor}\label{cor2}
Let $AN_n$ be the $n$-dimensional alternating group network for $n\geq 6$. Then $t_p(AN_n)=2n-5=\kappa_1(AN_n)$.
\end{cor}

\f\demo
%To prove the corollary, we only need to verify that $AN_n$ satisfies conditions in Theorem~\ref{main}.
Note that $N=|V(AN_n)|=\frac{n!}{2}\geq4(n-1)-2$ for $n\geq 6$, Condition $(1)$ in Theorem~\ref{main} holds.
Conditions $(2)$-$(4)$ in Theorem~\ref{main} hold by Lemmas~\ref{AN-1},~\ref{ANn4} and~\ref{ANn3}, respectively.
So $AN_n$ satisfies all conditions in Theorem~\ref{main}, and
$t_p(AN_n)=2n-5=\kappa_1(AN_n)$ for $n\geq 6$. \hfill\qed

\subsection{Application to BC networks}

Note that $L_{n}=\{X_{n}: X_{n} \ {\rm is\ an}\ n-{\rm dimensional\ BC\ network}\}$.
For a $BC$ network $X_n\in L_n$, the connectivity is $k=n\geq 5$, $l=0$, $N=|V|=2^n\geq 4n-2$ for $n\geq 5$ in Theorem~\ref{main}.
As a directive corollary of Theorem~\ref{main}, we can get the result $\kappa_1(X_n)=t_p(X_n)=2n-2$ in which
Zhu~\cite{Z08} determined $\kappa_1(X_n)=2n-2$ for $n\geq 4$. Fan and Lin~\cite{FX} obtained $t_p(X_n)=2n-2$ for $n\geq 4$.

\begin{lem}\label{bc3}
For any $X_n\in L_n$, if $U\subseteq V(X_n)$ with $2\leq|U|\leq 4n-8$ for $n\geq 3$, then $|N_{X_n}(U)|\geq 2n-2$.
\end{lem}

\f\demo We prove the lemma by using introduction on $n$.
If $n=3$, $2\leq|U|\leq 4n-8=4$, it is not difficult to see that $|N_{X_3}(U)|\geq 4$.
Assume that the lemma is true for $X_{m-1}$, where $m$ is an integer with $4\leq m\leq n-1$.
We consider $X_n$ for $n\geq 4$ as follows.

Since $X_n$ is $n$-regular $n$-connected triangle-free and $C(X_n)=2$, if $|U|=2$, then $|N_{X_n}(U)|\geq 2n-2$.
Now consider $3\leq |U|\leq 4n-8$. Note that $X_n$ contains two copies of $X_{n-1}$,
say $X_{n-1}^1$ and $X_{n-1}^2$, respectively. Let
$U_{i}=U\cap V(X_{n-1}^i)$ for $i\in\{1,2\}$.
Without loss of generality, we may assume that $|U_{1}|\geq |U_{2}|$.
It implies that $2\leq |U_{1}|$.

Case 1. $2\leq |U_{1}|\leq 4n-12$.
By the inductive hypothesis in $X_{n-1}^1$, $|N_{X_{n-1}^1}(U_{1})|\geq 2n-4$.
If $|U_{2}|=0$, then $U=U_1$.
$|N_{X_n}(U)|\geq |N_{X_{n-1}^1}(U_{1})|+|N_{\overline{X_{n-1}^1}}(U_{1})|\geq (2n-4)+2\geq 2n-2$.
If $|U_{2}|=1$, $|N_{X_{n-1}^2}(U_{2})|=\kappa(X_{n-1}^2)=n-1$.
Thus $|N_{X_n}(U)|\geq |N_{X_{n-1}^1}(U_{1})|+|N_{X_{n-1}^2}(U_{2})|\geq (2n-4)+(n-1)=3n-5\geq 2n-2$ for $n\geq 4$.
Now consider $2\leq |U_{2}|\leq |U_{1}|\leq 4n-12$ for $n\geq 4$, so $|N_{X_{n-1}^2}(U_{2})|\geq 2n-4$.
Thus, $|N_{X_n}(U)|\geq |N_{X_{n-1}^1}(U_{1})|+|N_{X_{n-1}^2}(U_{2})|\geq 2(2n-4)=4n-8 \geq 2n-2$ for $n\geq 4$.

Case 2. $4n-11\leq|U_{1}|\leq 4n-8$.

If $U=U_{1}$, by definition, $|N_{\overline{X_{n-1}^1}}(U_{1})|=|U_{1}|\geq 4n-11$.
Thus, $|N_{X_n}(U)|\geq |N_{\overline{X_{n-1}^1}}(U_{1})|\geq 4n-11\geq 2n-4$ for $n\geq 4$.
Now assume that $U\neq U_{1}$. Since the connectivity of $X_{n-1}^1$ is $n-1$ and
$|V(X_{n-1}^1)|-(4n-8)\geq \kappa(X_{n-1}^1)=n-1$ for $n\geq 4$, by Lemma~\ref{lem-3}, $|N_{X_{n-1}^1}(U_{1})|\geq n-1$.
Note that $U\neq U_{1}$ and $|U-U_{1}|\leq 3$, so $1\leq |U_{2}|\leq 3$.
If $|U_{2}|=1$, $|N_{X_{n-1}^2}(U_{2})|=\kappa(X_{n-1}^2)=n-1$.
Hence, $|N_{X_n}(U)|\geq |N_{X_{n-1}^1}(U_{1})|+|N_{X_{n-1}^2}(U_{2})|\geq 2n-2$ for $n\geq 4$.
Now suppose that $2\leq |U_{2}|\leq 3$. Since $|V(X_{n-1}^2)|-3\geq \kappa(B_{2})=n-1$ for $n\geq 4$,
by Lemma~\ref{lem-3}, $|N_{X_{n-1}^2}(U_{2})|\geq \kappa(X_{n-1}^2)=n-1$.
So $|N_{X_n}(U)|\geq |N_{X_{n-1}^1}(U_{1})|+|N_{X_{n-1}^2}(U_{2})|\geq 2n-2$ for $n\geq 4$.

By the above cases, the proof is completed.
\hfill\qed

\medskip
By Lemmas~\ref{bc2},~\ref{bc4} and~\ref{bc3} and Theorem~\ref{main}, we obtain the following Corollary~\ref{cor3}.

\begin{cor}\label{cor3}
For any $X_n\in L_n$, $t_p(X_n)=2n-2=\kappa_1(X_n)$ for $n\geq 5$.
\end{cor}

It is not difficult to check that the hypercube $Q_{n}$, the crossed cube $CQ_{n}$,
the M$\ddot{{\rm o}}$bius cubes $MQ_{n}$, the twisted cubes $TQ_{n}$ are all $n$-regular $n$-connected triangle-free BCs, then the following known result is derived directly.

\begin{cor}{\rm(\cite{FX})}\label{cor4}
Every pessimistic diagnosability of the hypercube $Q_{n}$, the crossed cube $CQ_{n}$, the
M$\ddot{o}$bius cubes $MQ_{n}$ and the twisted cubes $TQ_{n}$ is $2n-2$ for $n\geq 6$..
\end{cor}

\subsection{Application to the $k$-ary $n$-cube networks $Q_n^k$}

\begin{lem}\label{Qnk3}
Let $Q_{n}^{k}$ be a $k$-ary $n$-cube, where $k\geq2$ and $n\geq 1$ are integers.
\begin{enumerate}
\item [{\rm (1)}] For $n\geq 3$, let $U$ be a subset of $V(Q_n^2)$ with $2\leq|U|\leq4n-8$. Then $|N_{Q_{n}^{2}}(U)|\geq 2n-2$.
\item [{\rm (2)}] For $n\geq 3$, let $U$ be a subset of $V(Q_{n}^{3})$ and $2\leq |U|\leq 8n-10$, then $|N_{Q_{n}^{3}}(U)|\geq 4n-3$.
\item [{\rm (3)}] For $n\geq 3$ and $k\geq 4$, let $U$ be a subset of $ V(Q_{n}^{k})$ and $2\leq |U|\leq 8n-8$, then $|N_{Q_{n}^{k}}(U)|\geq 4n-2$.
\end{enumerate}
\end{lem}

\f\demo
Since the proof for the three cases are similar, we take $(2)$ as an example, the details for $(1)$ and $(3)$ are omitted.

Let $Q[0], Q[1], Q[2]$ represent the three disjoint subcubes obtained from $Q_{n}^{3}$ by partition over one dimension.
Let $U_{i}=U\cap V(Q[i])$ and $\overline{Q[i]}=Q_{n}^{3}-Q[i]$ for $i\in\{0,1,2\}$. Without loss of generality, we may assume that $|U_{0}|\geq |U_{1}|\geq |U_{2}|$.

The lemma is proved by the induction on $n$. When $n=3$, it is easy to check $|N_{Q_{3}^{3}}(U)|\geq 9$ for $2\leq |U|\leq 8n-10=14$.
We assume that the lemma is true for $Q_{m-1}^{3}$, where $m$ is an integer with $4\leq m\leq n-1$. We consider $Q_{n}^{3}$ for $n\geq 4$ as follows.

Case 1. $|U_{0}|\leq 1$.

In this case, $|U_{i}|\leq 1$ for all $0\leq i\leq 2$. Clearly, $2\leq |U|\leq 3$ because of $i\leq 2$. The Lemma follows if $|U|=2$ by Lemma~\ref{lem-k}.
Now assume that $|U|=3$. Since $Q_{n}^{3}$ is $2n$-regular and $Q[i]$ is isomorphic to $Q_{n-1}^{3}$, $|N_{Q_{n}^{3}}(U)|\geq 3\kappa(Q_{n-1}^{3})=3(2n-2)\geq 4n-3$ for $n\geq 3$.

Case 2. $2\leq |U_{0}|\leq 8n-18$.

By inductive hypothesis in $Q[0]$, $|N_{Q[0]}(U_{0})|\geq 4(n-1)-3=4n-7$.
If $U=U_0$, then $|N_{Q_{n}^{3}}(U)|=|N_{Q[0]}(U_{0})|+|N_{\overline{Q[0]}}(U_{0})| \geq 4n-7+2|U_0|\geq 4n-7+4=4n-3$. Assume $U\neq U_0$ in the following.
Note that $|U|\leq 8n-10$ and $|U_{0}|\geq |U_{1}|\geq |U_{2}|$, $|U_{1}|\leq 4n-5$.

If $|U_{1}|=1$, $|N_{Q[1]}(U_{1})|=\kappa(Q[1])=2n-2$.
Note that $Q[0]$ and $Q[1]$ are vertex disjoint, $|N_{Q_{n}^{3}}(U)|\geq |N_{Q[0]}(U_{0})|+|N_{Q[1]}(U_{1})|\geq (4n-7)+(2n-2)=6n-9\geq 4n-3$ for $n\geq 4$.
Now consider $2\leq |U_{1}|\leq 4n-5\leq 8n-18$ for $n\geq 4$, by inductive hypothesis in $Q[1]$, $|N_{Q[1]}(U_{1})|\geq 4(n-1)-3=4n-7$.
Thus, $|N_{Q_{n}^{3}}(U)|\geq |N_{Q[0]}(U_{0})|+|N_{Q[1]}(U_{1})|\geq 2(4n-7)=8n-14 \geq 4n-3$ for $n\geq 4$.

Case 3. $8n-17\leq|U_{0}|\leq 8n-10$.

If $U=U_{0}$, $|N_{Q_{n}^{3}}(U)|\geq |N_{\overline{Q[0]}}(U_{0})|=2|U_{0}|\geq 2(8n-17)\geq 4n-3$ for $n\geq 4$.
In the following, we assume the case of $U\neq U_{0}$. Since the connectivity of $Q[0]$ is $2n-2$, note that $U\neq U_0$, so $2\leq |U_0|\leq 8n-11$.
Since $|V(Q_{[0]})-U_0|=3^{n-1}-|U_0|\geq 3^{n-1}-(8n-11)\geq 2n-2=\kappa(Q[0])$ for $n\geq 4$, and by Lemma~\ref{lem-3}, $|N_{Q[0]}(U_{0})|\geq 2n-2$.

Note that $U\neq U_{0}$ and $|U-U_{0}|\leq 7$, so $1\leq |U_{1}|\leq 7$.

If $|U_{1}|=1$ and $|U_{2}|=0$, recall that the connectivity of $Q_{n}^{3}$ is $2n$ and $Q[i]$ is isomorphic to $Q_{n-1}^{k}$, $|N_{Q[1]}(U_{1})|=\kappa(Q[1])=2n-2$.
Note that each vertex in $Q[0]$ (resp. $Q[1]$) has an extra neighbor in $Q[2]$.
Hence, $|N_{Q_{n}^{3}}(U)|\geq |N_{Q[0]}(U_{0})|+|N_{Q[1]}(U_{1})|+|N_{Q[2]}(U_{0})|\geq 4n-4+(8n-17)=12n-21\geq 4n-3$ for $n\geq 4$.
If $|U_{i}|=1$ for $i=1,2$, $|N_{Q[i]}(U_{i})|=\kappa(Q[i])=2n-2$.
Hence, $|N_{Q_{n}^{3}}(U)|\geq |N_{Q[0]}(U_{0})|+|N_{Q[1]}(U_{1})|+|N_{Q[2]}(U_{2})|\geq 3(2n-2)=6n-6\geq 4n-3$ for $n\geq 4$.
Now suppose that $2\leq |U_{1}|\leq 7$.
Since $7<8n-17$ for $n\geq 4$, by inductive hypothesis in $Q[1]$, $|N_{Q[1]}(U_{1})|\geq 4(n-1)-3=4n-7$.
Thus, $|N_{Q_{n}^{3}}(U)|\geq |N_{Q[0]}(U_{0})|+|N_{Q[1]}(U_{1})|\geq (2n-2)+(4n-7)=6n-9\geq 4n-3$ for $n\geq 4$.

The proof is complete. \hfill\qed

\begin{rmk}\label{rmkQnk}
Esfahanian~\cite{E89} obtained $\kappa_1(Q_n^2)=2n-2$ for $n\geq 3$ and Day~\cite{D} got $\kappa_1(Q_n^3)=4n-3$, $\kappa_1(Q_n^k)=4n-2$ for $k\geq 4$. Kavianpour and Kim~\cite{KK} proved that $t_p(Q_n^2)=2n-2$ for $n\geq 3$ and  Wang et al.~\cite{WZF} derived
$t_p(Q_n^3)=4n-3$,  $t_p(Q_n^k)=4n-2$ for $k\geq 4$ and $n\geq 4$.
These results can be gotten directly as corollary of Theorem~\ref{main} as following.
\end{rmk}

Since $k^n\geq 4\kappa(Q_n^k)-2$ for $k\geq 3$ and $n\geq 3$ ($k=2$ and $n\geq 5$),
Condition $(1)$ in Theorem~\ref{main} holds. By Lemmas~\ref{Qnk2},~\ref{Qnk4} and~\ref{Qnk3}, Condition $(2)$-$(4)$ in Theorem~\ref{main} holds.

\begin{cor}\label{cor5}
Let $Q_{n}^{k}$ be a $k$-ary $n$-cube, where $k\geq2$ and $n\geq 1$ are integers. Then
\begin{enumerate}
\item [{\rm (1)}] $t_p(Q_n^2)=2n-2=\kappa_1(Q_n^2)$ for $n\geq 5$;
\item [{\rm (2)}] $t_p(Q_n^3)=4n-3=\kappa_1(Q_n^3)$ for $n\geq 3$;
\item [{\rm (3)}] $t_p(Q_n^k)=4n-2=\kappa_1(Q_n^k)$ for $n\geq 3$ and $k\geq 4$.
\end{enumerate}
\end{cor}

\subsection{Application to the split-star networks $S_{n}^{2}$}

Lin et al.~\cite{LX15} proved $\kappa_{1}(S_n^2)=4n-9$ for $n\geq 4$. However,
$t_{p}(S_n^2)$ has not been determined so far.
We can deduce the result by Theorem~\ref{main} in which for $S_n^2$, $k=2n-3$, $l=1$.

\begin{lem}\label{Sn23}
Let $S_{n}^{2}$ be the $n$-dimensional split-star network for $n\geq 4$.
If $U$ is a subset of $V(S_n^2)$ and $2\leq |U|\leq 8n-22$,
then $|N_{S_n^2}(U)|\geq 4n-9$.
\end{lem}

\f\demo
We prove the lemma by using the induction on $n$. Since $S_{4}^{2}$ is constructed by four disjoint triangles $S_{3}^2$, it is easy to verify that $|N_{S_{4}^{2}}(U)|\geq 7$ for $2\leq |U|\leq 10$.
By the inductive hypothesis, we assume that the lemma is true for $S_{m}^2$,
where $m$ is an integer with $5\leq m\leq n-1$.
Now we consider $S_n^2$.

Recall that $S_{n}^2$ is constructed by $n$ disjoint $S_{n-1}^2$s, denoted by $S_n^{2:i}$ for $i\in [n]$.
Let $U_{i}=U\cap V(S_n^{2:i})$ and $\overline{S_{n}^{2:i}}=S_{n}^2-S_{n}^{2:i}$ for $i\in [n]$.
Without loss of generality, we may assume that $|U_{1}|\geq |U_{2}|\geq \ldots \geq |U_{n}|$.
The following cases should be considered.

Case 1. $|U_{1}|\leq 1$.

In this case, $|U_{i}|\leq 1$ for all $i\in [n]$. Clearly, $2\leq |U|\leq n$ because of $U=\bigcup\limits_{i=1}^nU_i$.
If $|U|=2$, by Lemma~\ref{lem-k}, $|N_{S_n^2}(U)|\geq 2(2n-3)-2-1=4n-9$, the lemma follows.
Now assume that $3\leq |U|\leq n$. Since $S_{n}^2$ is $(2n-3)$-regular and $S_n^{2:i}$ is isomorphic to $S_{n-1}^2$, $|N_{S_n^2}(U)|\geq 3\kappa(S_n^{2:i})=3(2n-5)\geq 4n-9$ for $n\geq 5$.

Case 2. $2\leq |U_{1}|\leq 8n-30$.

By inductive hypothesis in $S_n^{2:1}$, $|N_{S_n^{1}}(U_{1})|\geq 4(n-1)-9=4n-13$.
Since $|U|\leq 8n-22$ and $|U_{1}|\geq |U_{2}|\geq \ldots \geq |U_{n}|$, $|U_{2}|\leq 4n-11$.
If $U=U_{1}$, by Lemma~\ref{Sn2}(2), $|N_{S_n^2}(U)|=|N_{S_n^{2:1}}(U_{1})|+|N_{\overline{S_n^{2:1}}}(U_{1})|\geq 4n-13+2|U_1|\geq 4n-9$.
Assume $U\neq U_1$ in the following. If $|U_{2}|=1$, $|N_{S_n^{2:1}}(U_{1})|=\kappa(S_n^{2:1})=2n-5$.
Note that $S_n^{2:1}$ and $S_n^{2:2}$ are vertex disjoint, $|N_{S_n^2}(U)|\geq |N_{S_n^{2:1}}(U_{1})|+|N_{S_n^{2:2}}(U_{2})|\geq 4n-13+2n-5=6n-18\geq 4n-9$ for $n\geq 5$.
Now consider $2\leq |U_{2}|\leq 4n-11$. Note that $4n-11\leq 8n-30$ for $n\geq 5$, by inductive hypothesis in $S_n^{2:2}$, $|N_{S_n^{2:2}}(U_{2})|\geq 4(n-1)-9=4n-13$.
Thus, $|N_{S_n^2}(U)|\geq |N_{S_n^{2:1}}(U_{1})|+|N_{S_n^{2:2}}(U_{2})|\geq 8n-26 \geq 4n-9$ for $n\geq 5$.

Case 3. $8n-29\leq|U_{1}|\leq 8n-22$.

By Lemma~\ref{Sn2}(2), $|N_{\overline{S_n^{2:1}}}(U_{1})|=2|U_{1}|$.
If $U=U_{1}$, $|N_{S_n^2}(U)|\geq|N_{\overline{S_n^{2:1}}}(U_{1})|=2|U_{1}|\geq 16n-58\geq 4n-9$ for $n\geq 5$.
In the following, we assume the case of $U\neq U_{1}$.
Since the connectivity of $S_n^{2:1}$ is $2n-5$,
and $(n-1)!-|U_{1}|\geq 2n-5=\kappa(S_n^{2:1})$ for $n\geq 5$,  by Lemma~\ref{lem-3}, $|N_{S_n^{2:1}}(U_{1})|\geq 2n-5$.
Note that $U\neq U_{1}$ and $|U-U_{1}|\leq 7$, so $1\leq |U_{2}|\leq 7$.

If $|U_{2}|=1$, recall that $S_n^2$ is $(2n-3)$-regular
and $S_n^{2:i}$ is isomorphic to $S_{n-1}^2$, $|N_{S_n^{2:2}}(U_{2})|=\kappa(S_n^{2:2})=2n-5$.
Hence, $|N_{S_n^2}(U)|\geq |N_{\overline{S_n^{2:1}}}(U_{1})|-|U-U_1|
\geq 16n-65\geq 4n-9$
for $n\geq 5$. Now suppose that $2\leq |U_{2}|\leq 7$.
Since $7\leq 8n-30$ for $n\geq 5$, by inductive hypothesis in $S_n^{2:1}$,
$|N_{S_n^{2:1}}(U_{1})|\geq 4(n-1)-9=4n-13$.
Thus, $|N_{S_n^2}(U)|\geq |N_{S_n^{2:1}}(U_{1})|+|N_{S_n^{2:2}}(U_{2})|\geq (2n-5)+(4n-13)=6n-18\geq 4n-9$ for $n\geq 5$.

By the above cases, the lemma holds. \hfill\qed

\begin{cor}\label{cor6}
Let $S_{n}^{2}$ be the $n$-dimensional split-star network for $n\geq 4$. Then
$t_{p}(S_n^2)=4n-9=\kappa_{1}(S_n^2)$.
\end{cor}

\f\demo
To prove the theorem, we only need to verify that $S_n^2$ satisfies conditions in Theorem~\ref{main}.
Note that $k=2n-3\geq 5$ for $n\geq 4$, $l=1$,
$N=|V(S_n^2)|=n!\geq4(2n-3)-2$ for $n\geq 4$, Condition $(1)$ in Theorem~\ref{main} holds.
By Lemmas~\ref{Sn2} and~\ref{Sn23}, Conditions $(2)$-$(3)$ in Theorem~\ref{main} holds. Condition $(4)$ holds by Lemma~\ref{Sn24}. $S_n^2$ satisfies all conditions in Theorem~\ref{main}, and thus
$t_p(S_n^2)=4n-9=\kappa_1(S_n^2)$.
\hfill\qed

\subsection{Application to the Cayley graphs generated by transposition trees $\Gamma_n$}

Let $\Gamma_n$ be Cayley graphs generated by transposition trees.
Yang et al.~\cite{Y10} determined $\kappa_1(\Gamma_n)=2n-4$ for $n\geq 3$.
However, $t_p(\Gamma_n)$ has not been known so far. By Theorem~\ref{main},
we immediately the following result which contains the above result.
Note that for $\Gamma_n$, $k=n-1$, $l=0$ in Theorem~\ref{main}.

\begin{lem}\label{Gamma3}
Let $\Gamma_n$ be Cayley graphs generated by transposition trees for $n\geq 4$.
If $U$ is a subset of $ V(\Gamma_n)$ and $2\leq |U|\leq 4n-12$, then $|N_{\Gamma_n}(U)|\geq 2n-4$.
\end{lem}

\f\demo
The lemma is proved by induction on $n$. When $n=4$, it is easy to check $|N_{\Gamma_n}(U)|\geq 4$
for $2\leq |U|\leq 4n-12=4$. We assume that the lemma is true for $\Gamma_{m}$, where $m$ is an integer with $4\leq m\leq n-1$. We consider $\Gamma_n$ for $n\geq 5$ as follows.

Recall that $\Gamma_n$ can be decomposed into $n$ copies of $\Gamma_{n-1}'s$, namely $\Gamma_{n-1}^{1}, \Gamma_{n-1}^{2}, \ldots,\Gamma_{n-1}^{n}$.
Let $U_{i}=U\cap V(\Gamma_{n-1}^{i})$ and $\overline{\Gamma_{n-1}^{i}}=\Gamma_n-\Gamma_{n-1}^{i}$ for $i\in[n]$. Without loss of generality, we may assume that $|U_{1}|\geq |U_{2}|\geq |U_{3}|\geq \ldots\geq |U_{n}|$.

Case 1. $|U_{1}|\leq 1$.

In this case, $|U_{i}|\leq 1$ for all $1\leq i\leq n$. Since $|U|\geq 2$, it implies $|U_1|=|U_2|=1$.
Since $\Gamma_n$ is $(n-1)$-regular and $\Gamma_{n-1}^{i}$ is isomorphic to $\Gamma_{n-1}$,
$|N_{\Gamma_n}(U)|\geq 2\kappa(\Gamma_{n-1}^{i})=2(n-2)=2n-4$ for $n\geq 5$.

Case 2. $2\leq |U_{1}|\leq 4n-16$.

By inductive hypothesis in $\Gamma_{n-1}^{1}$, $|N_{\Gamma_{n-1}^{1}}(U_{1})|\geq 2(n-1)-4=2n-6$.
Note that $|U_i|\leq |U_1|\leq 4n-16$ for $i\in \{2, 3, \ldots, n\}$.
If $|U_{2}|=1$, $|N_{\Gamma_{n-1}^{2}}(U_{2})|\geq \kappa(\Gamma_{n-1}^{2})=n-2$,
so $|N_{\Gamma_n}(U)|\geq |N_{\Gamma_{n-1}^{1}}(U_{1})|+|N_{\Gamma_{n-1}^{2}}(U_{2})|\geq (2n-6)+(n-2)=3n-8\geq 2n-4$ for $n\geq 5$.
If $2\leq |U_{2}|\leq 4n-16$, by inductive hypothesis in $\Gamma_{n-1}^{2}$, $|N_{\Gamma_{n-1}^{2}}(U_{2})|\geq 2(n-1)-4=2n-6$.
Thus, $|N_{\Gamma_n}(U)|\geq |N_{\Gamma_{n-1}^{1}}(U_{1})|+|N_{\Gamma_{n-1}^{2}}(U_{2})|\geq 2(2n-6)=4n-12 \geq 2n-4$ for $n\geq 5$.
Now consider $|U_{2}|=0$, then $|U_{i}|=0$ for $i\in \{3,4,\ldots, n\}$, it implies that $U=U_1$.
So $|N_{\Gamma_n}(U)|\geq |N_{\Gamma_{n-1}^{1}}(U_{1})|+|N_{\overline{\Gamma_{n-1}^{1}}}(U_1)|\geq 2n-6+|U_1|\geq 2n-6+2=2n-4$ for $n\geq 5$.

Case 3. $4n-15\leq|U_{1}|\leq 4n-12$.

If $U=U_{1}$, by Lemma~\ref{Gamma}, $|N_{\overline{\Gamma_{n-1}^{1}}}(U_{1})|=|U_{1}|\geq 4n-15$.
Since $(n-1)!-(4n-12)\geq n-2$ for $n\geq 5$, by Lemma~\ref{lem-3}, $|N_{\Gamma_{n-1}^{1}}(U_{1})|\geq \kappa(\Gamma_{n-1}^{1})=n-2$. Thus,
$|N_{\Gamma_n}(U)|=|N_{\overline{\Gamma_{n-1}^{1}}}(U_{1})|+|N_{\Gamma_{n-1}^{1}}(U_{1})|\geq 4n-15+(n-2)=5n-17\geq 2n-4$ for $n\geq 5$.
In the following, we assume that $U\neq U_{1}$.
It implies that $|U-U_{1}|\leq 3$, so $1\leq |U_{2}|\leq |U|-|U_1|\leq 3$.

If $|U_{2}|=1$, $|N_{\Gamma_{n-1}^{2}}(U_{2})|=\kappa(\Gamma_{n-1}^{2})=n-2$.
Recall that $|N_{\Gamma_{n-1}^{1}}(U_{1})|\geq n-2$.
Hence, $|N_{\Gamma_n}(U)|\geq |N_{\Gamma_{n-1}^{0}}(U_{0})|+|N_{\Gamma_{n-1}^{1}}(U_{1})|\geq 2n-4$ for $n\geq 5$.
Now suppose that $2\leq |U_{2}|\leq 3$. Since $|U_{2}|\leq 3\leq 4n-16$ for $n\geq 5$, by inductive hypothesis in $\Gamma_{n-1}^{2}$, $|N_{\Gamma_{n-1}^{2}}(U_{2})|\geq 2(n-1)-4=2n-6$.
Thus, $|N_{\Gamma_n}(U)|\geq |N_{\Gamma_{n-1}^{1}}(U_{1})|+|N_{\Gamma_{n-1}^{2}}(U_{2})|\geq (n-2)+(2n-6)=3n-8\geq 2n-4$ for $n\geq 5$.

By the above cases, the proof is completed.
\hfill\qed

\begin{cor}\label{cor-Tn}
Let $\Gamma_n$ be Cayley graphs generated by transposition trees for $n\geq 6$.
Then $t_p(\Gamma_n)=2n-4=\kappa_1(\Gamma_n)$ for $n\geq 6$.
\end{cor}

\f\demo
Note that $k=n-1\geq 5$ and $N=|V(\Gamma_n)|=n!\geq 4(n-1)-2$ for $n\geq 6$, Condition $(1)$ in Theorem~\ref{main} holds.
By Lemma~\ref{Gamma} and~\ref{Gamma3}, Condition $(2)$-$(3)$ in Theorem~\ref{main} holds.
Condition $(4)$ holds by Lemma~\ref{Gamma4}.
Thus, $\Gamma_n$ satisfies all conditions in Theorem~\ref{main}, $t_p(\Gamma_n)=2n-4=\kappa_1(\Gamma_n)$ for $n\geq 6$. \hfill\qed

\medskip

Since the star graph and the bubble-sort graph are Cayley graph generated by transposition trees,
The following corollary is gotten directly from Corollary~\ref{cor-Tn}.

\begin{cor}\label{cor8}
Let $S_n$ and $B_n$ are the star graph and the bubble sort graph,
then $t_p(S_n)=2n-4=\kappa_1(S_n)$ for $n\geq 6$, and $t_p(B_n)=2n-4=\kappa_1(B_n)$ for $n\geq 6$.
\end{cor}

\subsection{Application to the Cayley graphs generated by $2$-trees $\Gamma_{n}(\Delta)$}

\begin{lem}\label{tr3}
Let $\Gamma_{n}(\Delta)$ be a Cayley graph generated by the $2$-tree $\Delta$.
For $n\geq 4$, let $U$ be a subset of $V(\Gamma_{n}(\Delta))$ and $2\leq |U|\leq 8n-26$.
Then, $|N_{\Gamma_{n}(\Delta)}(U)|\geq 4n-11$.
\end{lem}

\f\demo The lemma is proved by the induction on $n$.
Since $\Gamma_{4}(\Delta)$ is constructed by $4$ disjoint triangles,
it is easy to verify that $|N_{\Gamma_{4}(\Delta)}(U)|\geq 5$ for $2\leq |U|\leq 7$.
By the inductive hypothesis, we assume that the lemma is true for $\Gamma_{m}(\Delta)$,
where $m$ is an integer with $5\leq m\leq n-1$.

Note that $\Gamma_{n}(\Delta)$ is constructed by $n$ disjoint $\Gamma_{n-1}(\Delta)$,
denoted by $\Gamma_n^{i}(\Delta)$ for $i\in [n]$.
Let $U_{i}=U\cap V(\Gamma_{n-1}^{i}(\Delta))$ and
$\overline{\Gamma_{n-1}^{i}(\Delta)}=\Gamma_{n}(\Delta)-\Gamma_{n-1}^{i}(\Delta)$ for $i\in [n]$.
Without loss of generality, we may assume that $|U_{1}|\geq |U_{2}|\geq \ldots \geq |U_{n}|$.
The following three cases should be considered.

Case 1. $|U_{1}|\leq 1$.

In this case, $|U_{i}|\leq 1$ for all $i\in [n]$. Clearly, $2\leq |U|\leq n$ because of $i\leq n$.
The Lemma follows if $|U|=2$ by Lemma~\ref{lem-k}.
Now assume that $3\leq |U|\leq n$. Since $\Gamma_{n}(\Delta)$ is $(2n-4)$-regular and $\Gamma_{n-1}^{i}(\Delta)$ is isomorphic to $\Gamma_{n-1}(\Delta)$,
$|N_{\Gamma_{n}(\Delta)}(U)|\geq 3\kappa(\Gamma_{n-1}^{i}(\Delta))=3(2n-6)\geq 4n-11$ for $n\geq 5$.

Case 2. $2\leq |U_{1}|\leq 8n-34$.

By inductive hypothesis in $\Gamma_{n-1}^{1}(\Delta)$, $|N_{\Gamma_{n-1}^{1}(\Delta)}(U_{1})|\geq 4(n-1)-11=4n-15$.
If $U=U_{1}$, $|N_{\Gamma_{n}(\Delta)}(U)|=|N_{\Gamma_{n-1}^{1}(\Delta)}(U_{1})|+
|N_{\overline{\Gamma_{n-1}^{1}(\Delta)}}(U_{1})|\geq 4n-15+2|U_1|\geq 4n-11$.
Assume $U\neq U_1$ in the following.
If $|U_{2}|=1$, $|N_{\Gamma_{n-1}^{2}(\Delta)}(U_{2})|=\kappa(\Gamma_{n-1}^{2}(\Delta))=2n-6$.
Note that $\Gamma_{n-1}^{1}(\Delta)$ and $\Gamma_{n-1}^{2}(\Delta)$ are vertex disjoint,
$|N_{\Gamma_{n}(\Delta)}(U)|\geq |N_{\Gamma_{n-1}^{1}(\Delta)}(U_{1})|+|N_{\Gamma_{n-1}^{2}(\Delta)}(U_{2})|
\geq 4n-15+(2n-6)\geq 6n-21$ for $n\geq 5$.
Now consider $2\leq |U_{2}|\leq |U_{1}|\leq 8n-34$,
by inductive hypothesis in $\Gamma_{n-1}^{2}(\Delta)$, $|N_{\Gamma_{n-1}^{2}(\Delta)}(U_{2})|\geq 4(n-1)-11=4n-15$.
Thus, $|N_{\Gamma_{n}(\Delta)}(U)|\geq |N_{\Gamma_{n-1}^{1}(\Delta)}(U_{1})|+
|N_{\Gamma_{n-1}^{2}(\Delta)}(U_{2})|\geq 8n-30\geq 4n-11$ for $n\geq 5$.

Case 3. $8n-33\leq|U_{1}|\leq 8n-26$.

By Lemma~\ref{prop-1}, $|N_{\overline{\Gamma_{n-1}^{1}(\Delta)}}(U_{1})|=2|U_{1}|$.
It is clear that the lemma holds if $U=U_{1}$.
In the following, we assume the case of $U\neq U_{1}$.
Since the connectivity of $\Gamma_{n-1}^{1}(\Delta)$ is $2n-6$, and
%and $\frac{(n-1)!}{2}-|U_{1}|\geq 2n-6=\kappa(\Gamma_{n-1}^{1}(\Delta))$ for $n\geq 5$,
by Lemma~\ref{lem-3},
$|N_{\Gamma_{n-1}^{1}(\Delta)}(U_{1})|\geq 2n-6$.
Note that $U\neq U_{1}$ and $|U-U_{1}|\leq 7$, so $1\leq |U_{2}|\leq 7$.

If $|U_{2}|=1$, % recall that $\Gamma_{n}(\Delta)$ is $(2n-4)$-regular and $\Gamma_{n-1}^{1}(\Delta)$ is
%isomorphic to $\Gamma_{n-1}(\Delta)$, $|N_{\Gamma_{n-1}^{2}(\Delta)}(U_{2})|=\kappa(\Gamma_{n-1}^{2}(\Delta))=2n-6$.
%Hence,
$|N_{\Gamma_{n}(\Delta)}(U)|\geq|N_{\Gamma_{n-1}^{1}(\Delta)}(U_{1})|+
|N_{\overline{\Gamma_{n-1}^{1}(\Delta)}}(U_{1})|-|U-U_1|\geq (2n-6)+2|U_1|-7\geq 18n-79\geq 4n-11$ for $n\geq 5$.
Now suppose that $2\leq |U_{2}|\leq 7$.
Since $7\leq 8n-32$ for $n\geq 5$, by inductive hypothesis in $\Gamma_{n-1}^{2}(\Delta)$,
$|N_{\Gamma_{n-1}^{2}(\Delta)}(U_{2})|\geq 4(n-1)-11=4n-15$.
Thus, $|N_{\Gamma_{n}(\Delta)}(U)|\geq |N_{\Gamma_{n-1}^{1}(\Delta)}(U_{1})|+|N_{\Gamma_{n-1}^{2}(\Delta)}(U_{2})|\geq (2n-6)+(4n-15)=6n-21\geq 2n-5$ for $n\geq 5$.

By the above cases, the lemma holds.\hfill\qed

\begin{cor}\label{cor9}
Let $G=\Gamma_{n}(\Delta)$ be a Cayley graph generated by the $2$-tree $\Delta$ for $n\geq 5$.
Then $\kappa_{1}(G)=4n-11=t_p(G)$.
\end{cor}

\f\demo
Note that $k=2n-4\geq 5$ and $\frac{n!}{2}\geq 4(2n-4)-2$ for $n\geq 5$, Condition $(1)$ in Theorem~\ref{main} holds.
By Lemma~\ref{prop-1} and~\ref{tr3}, Condition $(2)$ and $(3)$ in Theorem~\ref{main} holds.
Condition $(4)$ holds by $|F|\leq2k-3-l=2(2n-4)-3-1=4n-12<4n-11$ and Lemma~\ref{Delta}.
Thus, $\Gamma_{n}(\Delta)$ satisfies all conditions in Theorem~\ref{main}, and so
$t_p(\Gamma_{n}(\Delta))=4n-11=\kappa_1(\Gamma_{n}(\Delta))$ for $n\geq 5$. \hfill\qed

\subsection{Application to the burnt pancake networks $BP_n$}

\begin{lem}\label{BPn3}
Let $BP_n$ be the $n$-dimensional burnt pancake network.
For $n\geq 3$, let $U$ be a subset of $ V(BP_{n})$ and $2\leq |U|\leq 4n-8$, then $|N_{BP_{n}}(U)|\geq 2n-2$.
\end{lem}

\f\demo If $|U|=2$, by Lemma~\ref{BPn} and Lemma~\ref{lem-k},
for any two distinct vertices $u$ and $v$, so $|N_{BP_{n}}(U)|\geq 2n-2$.

Recall that $BP_{n}$ can be decomposed into $2n$ copies of $BP_{n-1}$'s, namely $BP_{n-1}^{i}$, for $i\in\langle n\rangle$.
Let $U_{i}=U\cap V(BP_{n-1}^{i})$ and $\overline{BP_{n-1}^{i}}=BP_n-BP_{n-1}^{i}$ for $i\in\langle n\rangle$.
Without loss of generality, we may assume that $|U_{1}|\geq |U_{2}|\geq |U_{3}|\geq \ldots\geq |U_{n}|\geq |U_{\bar{n}}|\geq |U_{\overline{n-1}}|\geq |U_{\bar{1}}|$.

The lemma is proved by using the induction on $n$.
If $n=3$, it is easy to check $|N_{BP_{n}}(U)|\geq 4$ for $2\leq |U|\leq 4n-8=4$.
We assume that the lemma is true for $BP_{m}$, where $m$ is an integer with $4\leq m\leq n-1$. We consider $BP_{n}$ for $n\geq 4$ as follows.

Case 1. $|U_{1}|\leq 1$.

In this case, $|U_{i}|\leq 1$ for all $1\leq i\leq n$. Since $|U|\geq 2$, it implies that $|U_1|=|U_2|=1$.
Since $BP_n$ is $n$-regular and $BP_{n-1}^{i}$ is isomorphic to $BP_{n-1}$, $|N_{BP_n}(U)|\geq 2\kappa(BP_{n-1}^{i})=2(n-1)= 2n-2$ for $n\geq 4$.

Case 2. $2\leq |U_{1}|\leq 4n-12$.

By inductive hypothesis in $BP_{n-1}^{1}$, $|N_{BP_{n-1}^{1}}(U_{1})|\geq 2(n-1)-2=2n-4$.
Note that $|U_i|\leq |U_1|\leq 4n-12$ for $i\in [n]\setminus \{1\}$.
If $U=U_{1}$, $|N_{BP_n}(U)|=|N_{BP_{n-1}^1}(U_1)|+|N_{\overline{BP_{n-1}^1}}(U_1)|\geq 4n-12+|U_1|\geq 4n-11$.
Assume $U\neq U_1$ in the following.
If $|U_{2}|=1$, $|N_{BP_{n-1}^{2}}(U_{2})|\geq \kappa(BP_{n-1}^{2})=n-1$,
so $|N_{BP_{n}}(U)|\geq |N_{BP_{n-1}^{1}}(U_{1})|+|N_{BP_{n-1}^{2}}(U_{2})|\geq (2n-4)+(n-1)=3n-5\geq 2n-2$ for $n\geq 4$.
If $2\leq |U_{2}|\leq 4n-12$, by inductive hypothesis in $BP_{n-1}^{2}$, $|N_{BP_{n-1}^{2}}(U_{2})|\geq 2(n-1)-2=2n-4$.
Thus, $|N_{BP_{n}}(U)|\geq |N_{BP_{n-1}^{1}}(U_{1})|+|N_{BP_{n-1}^{2}}(U_{2})|\geq 2(2n-4)=4n-8 \geq 2n-2$ for $n\geq 4$.

Case 3. $4n-11\leq|U_{1}|\leq 4n-8$.

Since $(n-1)!-(4n-8)\geq n-1$ for $n\geq 5$, by Lemma~\ref{lem-3},
$|N_{BP_{n-1}^{1}}(U_{1})|\geq \kappa(BP_{n-1}^{1})=n-1$.
If $U=U_{1}$, by Lemma~\ref{BPn}, $|N_{\overline{BP_{n-1}^{1}}}(U_{1})|=|U_{1}|\geq 4n-11$. Thus,
$|N_{BP_{n}}(U)|=|N_{\overline{BP_{n-1}^{1}}}(U_{1})|+|N_{BP_{n-1}^{1}}(U_{1})|\geq 4n-11+(n-1)=5n-2\geq 2n-2$
for $n\geq 4$. In the following, we assume that $U\neq U_{1}$.
It implies that $|U-U_{1}|\leq 3$, so $1\leq |U_{2}|\leq |U|-|U_1|\leq 3$.

If $|U_{2}|=1$, $|N_{BP_{n-1}^{2}}(U_{2})|=\kappa(BP_{n-1}^{2})=n-1$. Recall that $|N_{BP_{n-1}^{1}}(U_{1})|\geq n-1$.
Hence, $|N_{BP_{n}}(U)|\geq |N_{BP_{n-1}^{1}}(U_{1})|+|N_{BP_{n-1}^{2}}(U_{2})|\geq 2n-2$ for $n\geq 4$.
Now suppose that $2\leq |U_{2}|\leq 3$. Since $|U_{2}|\leq 3\leq 4n-12$ for $n\geq 4$, by inductive hypothesis in $BP_{n-1}^{2}$, $|N_{BP_{n-1}^{2}}(U_{2})|\geq 2(n-1)-2=2n-4$.
Thus, $|N_{BP_{n}}(U)|\geq |N_{BP_{n-1}^{1}}(U_{1})|+|N_{BP_{n-1}^{2}}(U_{2})|
\geq (n-1)+(2n-4)=3n-5\geq 2n-2$ for $n\geq 4$.

By the above cases, the proof is completed. \hfill\qed

\begin{rmk}
The extra connectivity of $BP_n$ was obtained by Song et al.~\cite{SXZC15}, $\kappa_1(BP_n)=2n-2$ for $n\geq 4$.
But $t_{p}(BP_n)$ is not known so far. By Theorem~\ref{main}, we immediately the following result which contains the above result.
\end{rmk}

\begin{cor}\label{cor10}
Let $BP_n$ be the $n$-dimensional burnt pancake network for $n\geq 5$. Then
$t_{p}(BP_n)=2n-2=\kappa_1(BP_n)$.
\end{cor}

\f\demo
Note that $k=n\geq 5$ and $N=|V(BP_n)|=n!\geq 4n-2$ for $n\geq 5$, Condition $(1)$ in Theorem~\ref{main} holds.
By Lemmas~\ref{BPn} and~\ref{BPn3}, Conditions $(2)$ and $(3)$ in Theorem~\ref{main} hold.
Condition $(4)$ holds by Lemma~\ref{BPn4}.
$BP_n$ satisfies all conditions in Theorem~\ref{main}, and so
$t_p(BP_n)=2n-2=\kappa_1(BP_n)$ for $n\geq 5$. \hfill\qed

\section{Concluding remarks}

This paper establishes the close relationship between these two parameter:
the extra connectivity and pessimistic diagnosability under the PMC model, by
proving $t_p(G)=\kappa_1(G)$ for some regular graphs $G$ with some conditions.
As applications, the pessimistic
diagnosability for each of split-star networks $S_n^2$, Cayley graphs generated by transposition trees $\Gamma_n$, Cayley graph generated by the $2$-tree $\Gamma_{n}(\Delta)$ and the burnt pancake networks $BP_n$ is gotten.
 As corollaries, the known results about the extra connectivity and the pessimistic diagnosability of many famous networks including the alternating group graphs~\cite{LZX},~\cite{T1}, the alternating group networks~\cite{Z99} , BC networks~\cite{Z08},~\cite{FX} and the $k$-ary $n$-cube networks~\cite{E89},~\cite{D},~\cite{KK},~\cite{WZF} are obtained directly.

%It is worth to investigate the relationship between the pessimistic diagnosability and the extra connectivity of a graph under the MM model.

\section*{Acknowledgments}

This work was supported by the National Natural Science Foundation of China (No.11371052, No.11271012 and No.11171020).

\end{document}